\documentstyle[12pt]{article}
\begin{document}
\newcommand{\n}{I{\hskip -5pt}N}
\newcommand{\r}{I{\hskip -5pt}R}
\newtheorem{definition}{Definition}[section]
\newtheorem{lemma}{Lemma}[section]
\newtheorem{theorem}{Theorem}[section]
\title{{The Existence of Global Solution for a Class of Semilinear Equations
 on Heisenberg Group
}
\thanks{This work is completed during the first author's visit to
the Institute of Theoretical Physics, Academia Sinica.}
\thanks{The second author is partially supported by the National Science Foundation
of P.R.China and Outstanding for Youth Foundation of Henan
Province and the National Science Foundation of Henan Province,
P.R.China} }
\author{{Zhujun Zheng}
\thanks{Email: zhengzj@henu.edu.cn }~~{Keping Lu}\\
{\scriptsize Mathematics Department of Henan University, Kaifeng,
475001}\\ {\scriptsize Institute of Mathematics, Henan University,
Kaifeng, 475001}\\ {Luo Xuebo}\\ {\scriptsize Institute of Applied
Mathematics, Northwestern Ploytechnical University, Xian, 710072}}
\date{}
\maketitle
\begin{center}
\begin{minipage}{120mm}
\vskip 1cm {\bf Abstract}{~Based on the concepts of a generalized
critical point and the corresponding generalized P.S. condition
introduced by Duong Minh Duc[1], we have proved a new $Z_2$ index
theorem and get a result on multiplicity of generalized critical
points. Using the result and  a quite standard variational method,
it is found that the equation $$ -\Delta_{H^n} u=|u|^{p-1} u
~~~~~~~~ x\in H^n $$ has infinite positive solutions. Our approach
can also be applied to study more general nonlinear problems.}

{\bf Key Words and Phrases~}Subelliptic Operator, $Z_2$ index,
Heisenberg Group, Generalized Critical Point.
\end{minipage}
\end{center}
\vskip 1cm
\baselineskip 20pt
\section{Introduction}
In this paper, we deal with the existence of the multiple global solutions
to the following nonlinear
equation
\begin{equation}
\label{1}
-\Delta_{H^n} u=|u|^{p-1} u
\end{equation}
where
$H^n$
is the Heisenberg group,
$
\Delta_{H^n}=\sum\limits_{i=1}^n(X_i^2 + Y_i^2)
$
is its subelliptic Laplacian operator. Under the real coordinate
$(x_1,\cdots,x_n,y_1,\cdots,y_n,t)$, the vector field $X_i$ and
$Y_i$ are defined by $$
\begin{array}{ll}
X_i=\frac{\partial}{\partial x_i}+2y_i\frac{\partial}{\partial t}&\\
&i=1,\cdots,n.\\
Y_i=\frac{\partial}{\partial x_i}-2x_i\frac{\partial}{\partial t}&\\
\end{array}
$$

 It is well known that $\{X_i, Y_i\}$ generate the real Lie
algebra of Lie group $H^n$ and  $$ [X_i,
Y_i]=4\delta_{ij}\frac{\partial}{\partial t}, i,j=1,\cdots,n. $$
In this Lie group, there is a group of natural dilations defined
by $$ \delta_\lambda(x,y,t)=(\lambda x, \lambda y, \lambda
t^2),\lambda
>0 $$ where $x=(x_1,\cdots,x_n), y=(y_1,\cdots,y_n)$. With this
group of dilations, the Lie group $H^n$
 is a two-step stratified nilpotent Lie group of homogeneous dimension $Q=2n+2$,
 and $\Delta_{H^n}$ is a homogeneous differential operator of degree 2.

Equation (1) comes from the CR-Yamabe problem(see [2]) and has
been studied by several authors(see [4], and the  references
therein). In  their works, they have got some results on the
existence of the boundary value problem of equation (1) on bounded
domain. In paper [2], when the domain is unbounded, they have
defined a thin condition, and proved a compact
Folland-Stein-Soblev type embedding theorem.
 By this compact theorem, they proved that the corresponding functional of equation (1) satisfies the P.S.
 condition, so by the
normal variational methods they gave some results about the boundary problem of equation (1).
 If the domain is unbounded and does not
satisfy the thin condition, to our knowledge, there exists no
report of progress on this problem up to now.

The global space $R^n$ is the most simple domain which does not
satisfy the thin condition. In the present paper, we study
equation (1) on the global space $R^n$, and find that it has
infinite positive solutions. The main result is the following



{\bf Main theorem} If $1<p<\frac{Q+2}{Q-2},$ where $Q=2n+2$, the
homogeneous dimension of $H^n$, then equation (1) has infinite
positive solutions belonging to $C^{2+\alpha}(H^n)$ with some
$\alpha > 0$.

Our idea is based on the concept of generalized critical point and
generalized P.S. condition introduced by   Duong Minh Duc[1]. In
his paper, he has obtained a deformation lemma and a generalized
Mountain Pass Lemma. Using these lemmas, he has studied a class of
nonlinear singular elliptic equations and obtained some existent
results. In our paper, based on the concept of generalized
critical point and generalized P.S. condition, we proved a new
deformation lemma. From our deformation lemma, we prove a
multi-existence result on a class of even functional which lacks
compact condition. Then using the lemma and quite standard
variational method, we present a proof of our main theorem.

It should be pointed out that our method can be used to study more
general problems, e.g. that related to a general unbounded domain
which has $Z_2$ symmetry and that with a more complicated
nonlinear term instead of $u|u|^{p-1}$ on the right hand of
equation (1).


\section{Preliminary results}
In this section, we give some definitions and lemmas related to the
generalized critical point, which come from  the ideal
of Duong Minh Duc(see [1]), and the definitions of $Z_2$ index. The theorem of $Z_2$
index related to generalized critical point is a new theorem and is our main tool in
the proof of main theorem.

Let $X$ be a Banach space. Assume that there exist a family of closed vector subspaces
$\{X_\rho\}_{\rho\in\cal{D}}$, and a family of linear maps $\{\pi_\rho\}_{\rho\in\cal{D}}$
such that $I={\rm Span}\{\bigcup\limits_{\rho\cal{D}}X_\rho\}$
 is dense in $X$, and $\pi_\rho(X)=X_\rho$.
\begin{definition}
\label{dy2.1}
Let $f$ be a continuous functional on $X$, $x\in X$.

(I). We call $x$ a generalized critical point of $f$ if there exists a sequence $\{x_j\}$ in $X$ such that
$$
\lim\limits_{j\rightarrow\infty}f'(x_j)=\lim\limits_{j\rightarrow\infty}\|\pi_\rho(x_j-x)\|=0
$$
for all $\rho\in\cal{D}$. In this case $\{x_j\}$ is called an approximation sequence of $x$.

(II). Let $x$ be a generalized critical point of $f$. We say that $x$ is regular if we can
find an approximation sequence $\{x_j\}$ of $x$ such that $\{f(x_j)\}$ is convergent in $\r$;
denote by $\bar{f}(x)$ the set of such limits. If $c$ in $\bar{f}(x)$, then $c$ is called a
generalized critical value of $f$.
\end{definition}
Denote by $K$ the set of all generalized critical points of $f$. For any real number $c$ and positive number $\varepsilon$, we put
$$
\begin{array}{l}
K_c=\{x\in K; c\in\bar{f}(x)\}\\
A(c,\varepsilon)=f^{-1}([c-\varepsilon, c+\varepsilon]).
\end{array}
$$ The generalized P.S. condition is defined as follows.
\begin{definition}
\label{dy2.2} Let $f$ be a continuous differential functional on
$X$. We say that $f$ satisfies generalized P.S. condition, if for
any sequence $\{x_j\}$ in $X$, along which $\{f(x_j)\}$ is bounded
and $\{f'(x_j)\}$ is convergent to 0, there exists a point $x\in
X$, and a subsequence $\{x_{j_k}\}$ of $\{x_j\}$, such that
$\lim\limits_{k\rightarrow\infty}\|\pi_\rho(x_{j_k}-x)\|=0$ for
all $\rho\in\cal{D}$.
\end{definition}
Next, we set $f$ to be a continuous differential functional, and $c\in\r$.
\begin{lemma}
\label{y2.1} Suppose that $K_c=\emptyset$. Then there exists
positive real number $b$ and $\varepsilon$, such that for any $x$
in $A(c,\varepsilon)$, we have $$ \|f'(x)\|>b. $$
\end{lemma}
Proof. If Lemma \ref{y2.1} is not true, there exist two sequences
of positive numbers $\{b_j\}, \{\varepsilon_j\}, b_j\rightarrow 0,
\varepsilon_j\rightarrow 0$, and $x_j\in A(c,\varepsilon_j)$, such
that $\|f'(x_j)\|<b_j$. Since $f$ satisfies generalized P.S.
condition, we have $c\in\bar{f}(x)$. That is $K_c\neq\emptyset$.
This contradicts to the condition $K_c=\emptyset$. So Lemma
\ref{y2.1} is true.
\begin{lemma}
\label{y2.2}
Suppose that $K_c\neq\emptyset$. Let $b$ and $\varepsilon $ be defined as
in Lemma \ref{y2.1}. Then there exist an subset $u$ of $X$ containing
$A(c,\varepsilon)$ and a locally Lipschitz continuous map $v$ from $u$
into $X$ such that for any $x$ in $u$
\begin{equation}
\label{2}
\|v(x)\|\leq 1
\end{equation}
\begin{equation}
\label{3}
f'(x) v(x)>\frac{1}{2}b.
\end{equation}
\end{lemma}
Proof. By Lemma \ref{y2.1}, for any $x$ in $A(c,\varepsilon)$, there
exists $h\in X, \|h\|=1$, such that
$
f'(x)h>\frac{1}{2}b.
$
Since $f'(x)$ is
continuous, there is a positive number $r_x$, such that for any $y\in
B(x,r_x)$, we have
\begin{equation}
\label{4}
f'(y)h>\frac{1}{2}b
\end{equation}
So we get a cover $B=\{B(x,r_x)|x\in A(c,\varepsilon)\}$ of
$A(c,\varepsilon)$. Let $\{B(x_i,r_{x_i})\}_{i\in Z}$ be a locally
finite subcover of $B$. Define the set $U=\bigcup\limits_{i\in
Z}B(x_i,r_{x_i})$ and the functional $q_i(x)$ which is the
distance from $x$ to $X\backslash B(x_i,r_{x_i})$. Then $q_i(x)$
is a Lipschitz continuous functional and $q_i|_{ X\backslash
B(x_i,r_{x_i})}\equiv 0$. Let $$
v(x)=\sum\limits_{i}\frac{q_i(x)}{\sum\limits_iq_j(x)}h_j $$ where
$h_j$ is defined by formula (\ref{4}). One can check that $v(x)$
is that we need. ~~~~~~~\#

Now, we can give our deformation lemma.
\begin{lemma}
\label{y2.3} Let $f$ be a continuous differential functional on
$X$, and $f$ satisfy the generalized P.S. condition. For $c\in\r,
K_c=\emptyset.$ $\varepsilon$ is the positive number given by
Lemma \ref{y2.1}.Then for every $q\in(0,\varepsilon)$, there is a
homomorphism $w$ on $X$, such that $$
\begin{array}{l}
(i) w(x)=x,~~~~~\forall x\in X\backslash A(c,\varepsilon),\\
(ii) w(f^{-1}((-\infty, c+q]))\subset f^{-1}((-\infty, c-q)).
 \end{array}
$$
\end{lemma}
Proof. Let $w$ to be a solution of Cauchy problem $$
\left\{\begin{array}{l} \frac{dw}{dt}=-v(w)\\ w(x,0)=x,
\end{array}
\right.
$$
where $v$ is defined as in Lemma 2.2. Then one can check that $w$ satisfies

(a). For any $s\in (0,t), w(x,s)\in A(c,q)$ implies
\begin{equation}
\label{5}
|w(x,t)-x|\leq t, f(x)-f(w(x,t))\geq \frac{1}{2}bt;
\end{equation}

(b). For any $x\in X\backslash A(c,\varepsilon), w(x,t)=x$;

(c). The function $f_1(t)=f(w(x,t))$ is not a increasing function on $t$.

Let $b$ be given as in Lemma \ref{y2.1}. It is clear that
$$
w(x,t)=x,{\hbox{for all}}~ x\in X\backslash A(c,\varepsilon)~\hbox{and}~
t\geq 0.
$$
Let $t_0$ be $\frac{4\varepsilon}{b}$. Observe that the trajectory $w(x,t)$
emanate from $x$ in $f^{-1}((-\infty, c+q))$, where $t\in [0,t_0]$.
If $x$ belongs to $f^{-1}((-\infty, -c-q)$, then $w(x,t)$ is in
$f^{-1}((-\infty, -c-\varepsilon))$.

Next it is proved that for any $x\in A(c,q), w(x,t_0)$ belongs to
$f^{-1}((-\infty, c-q])$. If this is not true, for any $s$ belongs
to $[0,t_0]$, one can find a point $x\in A(c,q)$. Then by the
formula (\ref{5}) we have $$
\begin{array}{ll}
f(w(x,0))-f(w(x,t_0))&=f(x)-f(w(x,t_0))\\
&\geq\frac{1}{2}bt_0=2\varepsilon>2q.
\end{array}
$$ This is a contradiction. Therefore the lemma follows.~~~~~~\#

When the set $K_c \neq \emptyset$, from the lemma 2.2 and the very
similar argument, we have the following deformation lemma.

\begin{lemma}
\label{yy2.4} Let $f$ be a continuous differential functional on
$X$, and $f$ satisfies the generalized P.S. condition. For $c\in
R$, $K_c \neq \emptyset$. Then there exists a positive
$\varepsilon$, and $q \in (0, \varepsilon)$, and a neighbourhood
$U$ of $K_c$,such that $f^{-1}(-\infty, c+q)\backslash U$ is a
deformation kernel of $f^{-1}(-\infty, c-\varepsilon)$.
\end{lemma}

Proof of lemma 2.4 is very similar to those of lemma 2.3 and the
corresponding lemma of deformation on the functional which
satisfies the P.S. condition. And when the functional has a $Z_2$
symmetry, the deformation can be chosen to be odd.

As the usual variational method of Ljusternik-Schnirelmann type
theory, after we have the above deformation lemma, we can give out
a $Z_2$ index on the generalized critical points. To be self
contained, we first give the definition of $Z_2$ index and some of
its properties we shall use. Then we give a theorem which is used
to compute the $Z_2$ index of the set of generalized critical
points.

Let $\cal{A}$ denote the set
$$
{\cal{A}}=\{A\in{\cal{A}}|A~{\hbox{is a symmetric close subset of }}X\}
$$
where symmetry means that $x\in A$ implies $-x\in A$. The
$Z_2$ index is defined as follows.
\begin{definition}
\label{dy2.3}
A function $i:{\cal{A}}\rightarrow Z_+\cup\{+\infty\}$ is called $Z_2$-index,
if for $A\in{\cal{A}}, i(A)$ is defined by

(I). If $A=\emptyset , i(A)=0$;

(II). if $A\neq\emptyset$, there exists a positive number $m$ and
a continuous odd map $\varphi :A\rightarrow\r^m\backslash\{0\}$,
then define $i(A)$ to be the minimum of this kind of $m$. i.e $$
i(A)=min\{m\in Z_+| {\hbox{there is a continuous odd map}}
\varphi:A\rightarrow\r^m\backslash\{0\}\} $$

(III). If $a\neq\emptyset$, and there is none positive integer satisfies (II),
define $i(A)=+\infty$.
\end{definition}
\begin{lemma}
\label{y2.4} The $Z_2$-index on ${\cal{A}}$ has the following
properties.

(I). $i(A)=0\Longleftrightarrow A=\emptyset$;

(II). If $A$ has only a pair of symmetric point, then $i(A)=1$;

(III). For any $A< B\in\cal{D}, A\subset B$, we have
$$
 i(A)\leq i(B);
$$

(IV). $i(A\cup B)\leq i(A)+i(B), ~~\forall A, B\in\cal{D}$;

(V). For any continuous odd map $\rho:X\rightarrow X$, and
$A\in{\cal{D}}, i(A)\leq i(\overline{\varphi(A)})$;

(VI). $A\in\cal{D}$, if $A$ is compact, then there exists a symmetric
neighborhood $N$, such that $i(\bar{N})=i(A)$, further more , if $A$ is compact,
then $i(A)<+\infty$;

(VII). Suppose $X_1$ is a subspace with dimension on $m, S$ is the
unit sphere, then $$ i(X_1\cap S)=m. $$
\end{lemma}

For their proofs, see e.g ref[7].

Let $M$ be a $C^{1,1}$ submanifold of $X, f$ is an even continuous
differential functional with generalized P.S. condition. If $k\leq
i(M)\leq \infty$, by Lemma \ref{y2.4}, the set $$
Z_k=\{A\in{\cal{D}} : A\subset M, i(A)\geq k\} $$ is not empty and
invariant under odd continuous map.

For any positive integer $k\leq r(M)$, define
\begin{equation}
\label{6}
C_k=\inf\limits_{A\in Z_k}\sup\limits_{u\in A} f(u).
\end{equation}
Then we have the following theorem.
\begin{theorem}
\label{d2.1}
Suppose $k<r(M), k+l-1\leq r(M)$. If
$$
-\infty<C_k=C_{k+1}=\cdots=C_{k+l-1}=c<+\infty
$$
then $i(K_c)\geq l$.
\end{theorem}
Proof. $K_c$ is a compact subset of $M$ (In fact, for any
$\{x_j\}\subset K_c$, since $f$ satisfies generalized P.S.
condition, there exist $x\in X$ and a subsequence
$\{x_{j_k}\}\subset\{x_j\}$, such that, for any $\rho\in\cal D$,
$$ \lim\limits_{k\rightarrow\infty}\|\pi_\rho(x_{j_k}-x\|=0. $$
Notice that the subspace ${\rm
Span}\bigcup\limits_{\rho\in{\cal{D}}} X_\rho$ of $X$ is dense in
$X$, so $\{x_{j_k}\}$ have a limit $x\in X$. From the continuity
of the functional $f$, we have $f(x)=c$ and $f'(x)=0$. Obviously
$x\in M$. Hence $x\in K_c$), so $i(K_c)$ can be defined. By Lemma
\ref{y2.4}, there exists a symmetric neighborhood $N$ in $M$, such
that $i(N)=i(K_c)$. If $i(K_c)<l$, by the deformation Lemma
\ref{yy2.4}, there is a positive number $\varepsilon>0$ and a
homomorphism $ w : X\rightarrow X$, such that $$
w(f_{c+\varepsilon}\backslash N)\subset f_{c-\varepsilon} $$ where
$f_{c+\varepsilon}$ is defined by $f_{c+\varepsilon}=\{x\in
M|f(x)<c+\varepsilon\}$, $f_{c-\varepsilon}$ is defined by the
same way. For the positive number $\varepsilon $ defined above, by
the definition $c_{m+l-1}$, we can find a set
$A\in{\cal{A}}_{m+l-1}$, such that
\begin{equation}
\label{7}
c\leq \sup f(x)<c+\varepsilon.
\end{equation}
So we have $$ \sup\limits_{x\in w(A\backslash
N)}f(x)<c-\varepsilon. $$ From Lemma \ref{y2.4}, noting $A\in
{\cal{A}}_{k+l-1}$, $$
\begin{array}{ll}
k+l-1&\leq i(A\backslash N)+i(\bar N)\\
&\leq i(w(A\backslash N)+l-1.
\end{array}
$$ So $$ r(A\backslash N)\geq k. $$
This implies $ A\backslash
N\subset\cal{A}$, which means
$$ \sup\limits_{x\in w(A\backslash
N)}f(x)\geq c_m=c. $$ This conflict with formula (\ref{6}). So
$r(K_c)\geq l$.~~~\#

As a result, we have the following theorem.
\begin{theorem}
\label{d2.2} Suppose $I$ is a continuous functional on $M.~ M$ is
a $C^{1,1}$ symmetric submanifold of $ X\backslash \{0\}$. Suppose
further that $I$ satisfies the generalized P.S. condition and
bounded from below. Then $I$ has at least $\widetilde{r}(M)$ pairs
of generalized critical points, where $\widetilde{i}(M)$ is
defined by
 $$
\widetilde{i}(M)=max\{i(A)|A\subset M, A~\hbox{is symmetric}\}.
$$
\end{theorem}
\section{Main theorem and its proof}
In this section, we present out main theorem and its proof. First
we define some space, some of which are not new. Then we give our
main theorem and its proof.

For any domain $\Omega\subset H^n$, let $C_0^\infty(\Omega)$ denote the
set of smooth functions with compact support. For $u\in C_0^\infty(\Omega)$,
define
\begin{equation}
\label{8}
\|u\|^2=\int_\Omega|\nabla_{H^n}u(x)|^2dx.
\end{equation}
One can check that $\|\cdot \|$ is a norm (see \cite{3}). Let $S_1^2(\Omega)$
denote the complete of $C_0^\infty(\Omega)$ under the norm $\|\cdot\|$.
Then we have the following embedding theorem \cite{3}.
\begin{lemma}
\label{y3.1}
Let $D\subset H^n$ be a bounded open set, then the embedding
$$
\stackrel{\circ}{S}_1^2(D)\hookrightarrow L^p(D), ~~~1\leq p<\frac{2Q}{Q-2}
$$
is compact.
\end{lemma}

In $H^n$, we define metric
$$
d((x,y,t),(x', y', c'))=[\bar t^2+(|\bar x|^2+|\bar y|^2)^2)^{\frac{1}{2}}
$$
where $(\bar x, \bar y, \bar t)=(x, y, t)-(x', y', t')$.
Let $B_k=\{x\in H^n|d((x,y,t),0)<k\}$. For $k>1$, define
$\rho_k\in C_0^\infty(B_{k+1})$, and $\rho_k|_{B_k}\equiv 1$. For $k>1$,
denote $H_k=\stackrel{\circ}{S}_1^2(B_k)$, and define a continuous map
$\pi_k : \stackrel{\circ}{S}_1^2(\r^n)\rightarrow H_k$ by
$$
\pi_k(u)=\rho_k u.
$$
In $\stackrel{\circ}{S}_1^2(\r^n)$, we define
$$
M=\{u\in \stackrel{\circ}{S}_1^2(\r^n)|\frac{1}{p+1}\int_{\r^n}|u|^{p+1}dx=1\}.
$$
It is obvious that $M$ is a $C^{1,1}$  submanifold of
$\stackrel{\circ}{S}_1^2(\r^n)$, and is symmetric , and that $0\not\in M$. We
define a functional on $M$
\begin{equation}
\label{9}
I(u)=\frac{1}{2}\int_{\r^n}|\nabla_{H^n}u|^2dx,~~~\forall u\in M.
\end{equation}
Then we have the following lemma.
\begin{lemma}
\label{y3.2} If $1\leq p<\frac{Q+2}{Q-2}$, the functional $I$
defined by formula (\ref{9}) satisfies generalized P.S. condition.
\end{lemma}
Proof. Obviously, $I$ is a continuous differential functional defined on $M$.
Along $\{u_j\}\subset M$, $\{I(u)\}$ is bounded and $\{I'(u)\}$
converges to 0.

$\forall k\in Z^+$, and $k\geq 2$, then $C_0^\infty(B_k)$ is dense in
$\stackrel{\circ}{S}_1^2(B_k)$, and $\{I(\pi_k(u_j))\}$ is bounded.
For any $u\in C_0^\infty(B_k)$, we have
$$
\begin{array}{ll}
<I'(\pi_k(u_j)),u>&=\int_{B_k}\nabla_{H^n}\pi_k(u_j)\cdot \nabla_{H^n} udx\\
\\
&=\int_{\r^n}\nabla_{H^n}(\rho_ku_j)\cdot\nabla_{H^n} udx\\
\\
&=\int_{\r^n}\rho_k\nabla_{H^n}u_j\cdot\nabla_{H^n}udx+\int_{\r^n}u_j
\nabla\rho_k\cdot\nabla_{H^n}udx\\
\\
&=\int_{B_k}\rho_k\nabla_{H^n}u_j\cdot\nabla_{H^n}udx+\int_{B_k}u_j
\nabla\rho_k\cdot\nabla_{H^n} udx\\
\\
&\leq C \int_{B_k}\nabla_{H^n}u_j\cdot\nabla_{H^n}udx\rightarrow 0,
\end{array}
$$ where $C$ is a positive constant independent of $k$. So on
$S_1^2(B_k), I^1(\pi_k(u_j))\rightarrow 0$.

This implies that, $\forall \varepsilon>0, \exists N$, if $r, s>N$,
we have, for any $k\in\cal{D}$,
\begin{equation}
\label{10}
\int\nabla_{H^n}\pi_k(u_r) \nabla_{H^n}(u_r-u_s)dx\leq\varepsilon\|
\pi_k(u_r-u_s)\|,
\end{equation}
\begin{equation}
\label{11}
\int\nabla_{H^n}\pi_k(u_s) \nabla_{H^n}
\pi_k(u_s-u_r)dx\leq\varepsilon\|\pi_k(u_r-u_s)\|.
\end{equation}
By the inequalities (\ref{10}),(\ref{11}), we have
$$
\|\pi_k(u_r-u_s)\|<\varepsilon ~~~~\forall k\in\cal D.
$$
So $\{\pi_k(u_i)\}$ converges in $H_k$ and converges almost every where in $B_k$.
But for $k+1, \pi_{k+1}(u_i)|_{B_k}=\pi_k(u_i)$. So we can choose a
function  $u\in S_1^2(\r^n)$, such that
\begin{equation}
\label{12}
\pi_k(u_i)\rightarrow\pi_k(u)~~~~~\forall k\in\cal D.
\end{equation}
The formula (\ref{12}) implies that the functional satisfies
generalized P.S. condition.

By Theorem \ref{d2.2} and Lemma \ref{y3.2}, we have the following theorem
\begin{theorem}
\label{d3.1}
If $1\leq p<\frac{Q+2}{Q-2}$, the functional (9) has infinite numbers of
generalized critical point.
\end{theorem}
Next, we shall prove that a generalized critical point of $I$ on $M$
is its critical point in $M$. This is our following lemma.
\begin{lemma}
\label{y3.3}
If $u$ is a generalized critical point of the functional $I$ on $M$,
it must be a critical point of $I$ on $M$.
\end{lemma}
Proof. Suppose that $u$ is a generalized critical point of $I$ on $M$.
By the definition of generalized critical point, there exists a sequence
$\{u_j\}\subset M$, such that, for any positive integer number $k$,
$$
I'(\pi_k(u_j))\rightarrow 0, \|\pi_k(u_j)-\pi_k(u)\|\rightarrow 0.
$$
Since $\bigcup\limits_{k=1}^\infty\pi_k(S^2_1(\r^n))\cap M$ is a dense
subset of $M$ and
$$
\|u_j-u\|\leq\|\pi_k(u_j)-\pi_k(u_j)\|+\|\pi_k(u_j)-\pi_k(u)\|+\|\pi_k(u)-u\|
$$
we have
$$
u_j\rightarrow u.
$$
For $I$ is a continuous differential functional, we have,
as $j\rightarrow \infty$,
$$
I'(\pi_k(u_j))\rightarrow I'(u).
$$
So we have $I'(u)=0$. This implies $u$ to be a critical point of $I$ on $M$.

By Lemma \ref{y3.3} and Lemma \ref{y3.2} we know that the functional $I$ has
infinitely many pairs of critical point on $M$; we denote them by
$\{\widetilde{u}_j\}$. So by Lagrangian multiplier rule,
there exists parameter $\mu\in\r$, such that
 $$
\begin{array}{ll}
<v, (DI(\widetilde{u}_j)-\mu\widetilde{u_j}|\widetilde{u_j}|^{p-2})>
&=\frac{1}{2}\int_{\r^n}\nabla_h\widetilde{u}_j\cdot\nabla_hv-
\mu\widetilde{u}_j|\widetilde{u}_j|^{p-2}vdx\\
&=0~ \hbox{for all }v\in S_1^2(H^n).
\end{array}
$$
Setting $v=\widetilde{u}_j$ in this equation yields that
$$
 \frac{1}{2}\int_{H^n}|\nabla_h\widetilde{u}_j|^2dx
 =\mu\int_{H^n}|\widetilde{u}_j|^{p+1}dx=(p+1)\mu.
$$
For $\widetilde{u}_j\neq 0$ in $S_1^2(H^n)$,
we have $\mu>0$. Rescaling $\widetilde{u_j}$ with a suitable power of $(p+1)\mu$,
we obtain a weak solution $u_j=((p+1)\mu)^{\frac{1}{p-2}}
\widetilde{u}_j\in\stackrel{\cdot}{S}_1^2(H^n)$ of problem
(\ref{1}), in the sense of
$$
\int_{H^n}
(\nabla_Hu_j\cdot\nabla v-|u_j|^{p-2}u_jvdx=0
~~\hbox{for all}
v\in\stackrel{\cdot}{S}_1^2(H^n).
$$
By the regulation theorem of operator $\Delta_{H^n}$(see [1], [2]),
we know $u_j\in C^{2+\alpha}(\r^n)$. Hence they are strong solutions
of problem (\ref{1}). Further more, by the minimum theorem of $\Delta_{H^n}$,
we know $u_j(p)\neq 0$ for all $p\in H^n$. So we get our main theorem.
\begin{theorem}
\label{d3.2}
If $1<p<\frac{Q+2}{Q-2}$, where $Q=2n+2$ is the homogeneous dimension of $H^n$, the problem
$$
\left\{\begin{array}{l}
-\Delta_{H^n}u=|u|^{p-1}u ~~x\in\r^n\\
u(x)>0
\end{array}
\right. $$ has infinite $C^{2+\alpha}$ solutions.
\end{theorem}

\end{document}